\documentclass{amsart}

\usepackage[english]{babel}
\usepackage[latin1]{inputenc}
\usepackage{bm}
\usepackage{fullpage}
\usepackage{amssymb}
\usepackage{amsmath}
\usepackage{latexsym}
\usepackage{graphicx}

\newtheorem{theorem}{Theorem}

\newtheorem{lemma}{Lemma}
\newtheorem{corollary}{Corollary}

\def\R{\mathbb{R}}
\def\Rt{\mathbb{R}_{\mbox{\textup{\tiny{max}}}}}

\author{Stéphane Gaubert}
\address{Inria Rocquencourt, Domaine de Voluceau, 78153 Le Chesnay,
France}
\email{stephane.gaubert@inria.fr}
\author{Frédéric Meunier}
\address{Université Paris Est, LVMT, ENPC, 6-8 avenue Blaise Pascal, Cité Descartes
Champs-sur-Marne, 77455 Marne-la-Vallée cedex 2, France.}
\email{frederic.meunier@enpc.fr}

\title{Carathéodory, Helly and the others in the max-plus world}

\begin{document}

\begin{abstract}
Carathéodory's,
Helly's and Radon's theorems are three basic results in discrete geometry.
Their max-plus counterparts have been proved by various authors.
In this paper, more advanced results in discrete geometry are shown to have also their max-plus counterparts:
namely,
the colorful Carathéodory theorem and the Tverberg theorem. A
conjecture connected to the Tverberg theorem -- Sierksma's
conjecture --, although still open for the usual convexity, is shown
to be true in the max-plus settings.
\end{abstract}

\maketitle

\section{Introduction}

Three basic theorems gave rise to this new topic that is discrete
geometry of convex sets, namely, Carathéodory's theorem, Helly's
theorem and Radon's theorem. We state them here for sake of
completeness.

\begin{theorem}[Carathéodory's theorem] \label{thm:cara} Suppose given $n\geq d+1$ points $\boldsymbol{x}_1,\boldsymbol{x}_2,\ldots,\boldsymbol{x}_n$ in $\R^d$ and a point $\boldsymbol{p}$ in $\mbox{\textup{conv}}\{\boldsymbol{x}_1,\boldsymbol{x}_2,\ldots,\boldsymbol{x}_n\}$. Then there is a subset $I\subseteq\{1,\ldots,n\}$ of cardinality $d+1$ such that $\boldsymbol{p}$ is in the convex hull of $\bigcup_{i\in I}\{\boldsymbol{x}_i\}$.
\end{theorem}

\begin{theorem}[Radon's theorem]
Let $X$ be a set of $d+2$ points in $\mathbb{R}^d$. Then there are
two pairwise disjoint subsets $X_1$ and $X_2$ of $X$ whose convex
hulls have a common point.
\end{theorem}

\begin{theorem}[Helly's theorem]
Let $\mathcal{F}$ be a finite collection of convex sets in $\R^d$. If every $d+1$ members of $\mathcal{F}$ have a nonempty intersection, then the whole collection have a nonempty intersection.
\end{theorem}

Recently, different authors have shown that max-plus versions of
these theorems exist. For Carathéodory's theorem, it is
explicitly stated in the paper of Develin and Sturmfels
\cite{DevelinSturmfels04}. A max-plus Helly theorem was proved in
the paper of Gaubert and Sergeev \cite{GaubertSergeev07}. When no
component of the vectors can be equal to $-\infty$, the max-plus
Radon theorem is then a consequence of Gondran-Minoux's theorem
\cite{GondranMinoux78,Butkovic03}. In the paper where Briec and Horvath introduce the notion of $\mathbb{B}$-convexity \cite{BriecHorvath04}, these three max-plus theorems are proved, too (without $-\infty$ components).

Other theorems have followed which of Carathéodory, Helly and
Radon,and the purpose of this paper is precisely to show that these
theorems have also max-plus versions. Section \ref{sec:colorCara} is
devoted to the max-plus version of the beautiful colorful
Carathéodory theorem proved by B\'ar\'any \cite{Barany82}.

\begin{theorem}[Colorful Carathéodory's theorem \cite{Barany82}] \label{thm:colorcara} Suppose given $d+1$ finite point sets $X_1, X_2, \ldots, X_{d+1}$ and a point $\boldsymbol{p}$ in $\R^d$ such that the convex hull of each $X_i$ contains $\boldsymbol{p}$, then there are $d+1$ points $\boldsymbol{x}_1,\boldsymbol{x}_2,\ldots,\boldsymbol{x}_{d+1}$ such that $\boldsymbol{x}_i\in X_i$ for each $i$ and such that the point $\boldsymbol{p}$ is the convex hull of the points $\boldsymbol{x}_1,\boldsymbol{x}_2,\ldots,\boldsymbol{x}_{d+1}$.
\end{theorem}

As said, the previous proofs of Radon's theorem do not deal with
$-\infty$ components; this is settled in Section \ref{sec:Radon}.
This leads to a new proof of the max-plus Helly theorem, presented
in Section \ref{sec:Helly}, since the usual way for proving the
Helly theorem is with the use of the Radon theorem.



Radon's theorem has a a beautiful generalization, Tverberg's
theorem. The max-plus version is proved in Section
\ref{sec:Tverberg}.

\begin{theorem}[Tverberg's theorem \cite{Tverberg66}]
Let $X$ be a set of $(d+1)(q-1)+1$ points in $\mathbb{R}^d$. Then
there are $q$ pairwise disjoint subsets $X_1, X_2,\ldots,X_q$ of $X$
whose convex hulls have a common point.
\end{theorem}

The case $q=2$ reduces to the usual Radon theorem.

A natural question is about the number of these partitions into $q$
subsets (each of these partition is called a {\em Tverberg
partition}). A famous conjecture is the following one, also called
the Dutch cheese conjecture, since Sierksma has offered a Dutch
cheese for a solution of this problem.

\medskip

\paragraph{\bf Conjecture} (Sierksma's conjecture) Let $q\geq 2$, $d\geq 1$ and put $N=(d+1)(q-1)$. For every $N+1$ points in $\R^d$ the number of unordered Tverberg partitions is at least $((q-1)!)^d$.

\medskip

This conjecture is still open. One can naturally ask whether this
conjecture holds in the max-plus settings. Surprisingly, it is
possible to prove it in this case, with a quite simple proof. This
is done in the last section of the paper.

\bigskip

\paragraph{\bf Notation: } Before starting, we introduce some notations. Note $\R\cup\{-\infty\}$ by $\Rt$. All sets, collections,... will be denote by capital letters, vectors with bold symbols (for example $\boldsymbol{x}$), and scalars with the usual typography.

If $\lambda$ is in $\Rt$ and
$\boldsymbol{x}=\left(\begin{array}{c}x_1 \\ \vdots \\
x_d\end{array}\right)$, then $\lambda+\boldsymbol{x}$ will denote
$\left(\begin{array}{c}\lambda+x_1 \\ \vdots \\
\lambda+x_d\end{array}\right)$.

The usual convex hull of the points
$\boldsymbol{x}_1,\ldots,\boldsymbol{x}_n$ in $\R^d$ is denoted by
$\mbox{\textup{conv}}\{\boldsymbol{x}_1,\ldots,\boldsymbol{x}_n\}$
and the max-plus convex hull of the points
$\boldsymbol{x}_1,\ldots,\boldsymbol{x}_n$ in $\Rt^d$ by
$\mbox{\textup{mpconv}}\{\boldsymbol{x}_1,\ldots,\boldsymbol{x}_n\}$.
This last set is then the set of points $\boldsymbol{x}$ in $\Rt^d$
such that there exist $\lambda_1,\ldots,\lambda_n$ in $\Rt$ such
that $\max_{i=1,\ldots,n}\lambda_i=0$ and
$\boldsymbol{x}=\max_{i=1,\ldots,n}(\lambda_i+\boldsymbol{x}_i)$
(for each of the $d$ components, one takes the maximum of the $n$
possible distinct values).

\section{The colorful Carathéodory theorems}

\label{sec:colorCara}

Before stating and proving the max-plus counterpart of the colorful
Carathéodory theorem, one sates an equivalent version,
straightforwardly derived from the colorful Carathéodory theorem.
Surprisingly, although the max-plus colorful Carathéodory theorem
has a very simple proof (like for the simple Carathéodory theorem),
the version with the convex set $C$ instead of the point
$\boldsymbol{p}$ needs in the max-plus settings more advanced tools.

This equivalent version is:

\begin{theorem}
Suppose given $d+1$ finite point sets $X_1, X_2, \ldots, X_{d+1}$
and a convex set $C$ in $\R^d$ such that the convex hull of each
$X_i$ intersects $C$. Then there are $d+1$ points
$\boldsymbol{x}_1,\boldsymbol{x}_2,\ldots,\boldsymbol{x}_{d+1}$ such
that $\boldsymbol{x}_i\in X_i$ for each $i$ and such that
$\mbox{\textup{conv}}\{\boldsymbol{x}_1,\boldsymbol{x}_2,\ldots,\boldsymbol{x}_{d+1}\}$
intersects the convex $C$.
\end{theorem}

\begin{proof}
For $i\in [d+1]$, one has $\boldsymbol{c}_i\in C$ such that
$\boldsymbol{c}_i\in \mbox{conv}(X_i)$. Hence one has
$\boldsymbol{0}\in\mbox{conv}\left(X'_i\right)$ where
$X'_i:=\{\boldsymbol{x}-\boldsymbol{c}_i:\,\boldsymbol{x}\in X_i\}$.
Applying Theorem \ref{thm:colorcara} to the sets
$X'_1,\ldots,X'_{d+1}$ and the point
$\boldsymbol{p}:=\boldsymbol{0}$ leads to the existence of $d+1$
points
$\boldsymbol{x}_1,\boldsymbol{x}_2,\ldots,\boldsymbol{x}_{d+1}$ such
that $\boldsymbol{x}_i\in X_i$ for each $i$ and such that
$\boldsymbol{0}\in\mbox{conv}\{\boldsymbol{x}_1-\boldsymbol{c}_1,\ldots,\boldsymbol{x}_{d+1}-\boldsymbol{c}_{d+1}\}$.
It implies that
$\mbox{conv}\{\boldsymbol{c}_1,\ldots,\boldsymbol{c}_{d+1}\}\cap\mbox{conv}\{\boldsymbol{x}_1,\ldots,\boldsymbol{x}_{d+1}\}\neq\emptyset$,
and hence that
$\mbox{conv}\{\boldsymbol{x}_1,\boldsymbol{x}_2,\ldots,\boldsymbol{x}_{d+1}\}$
intersects the convex $C$.
\end{proof}

We will see the max-plus version of this theorem later. Let us now
state the max-plus colorful Carathéodory theorem and prove it.

\begin{theorem}[Max-plus colorful Carathéodory's theorem]\label{thm:trop_colorcara}
Suppose given $d+1$ finite point sets $X_1, X_2, \ldots, X_{d+1}$
and a point $\boldsymbol{p}$ in $\Rt^d$ such that the max-plus
convex hull of each $X_i$ contains $\boldsymbol{p}$. Then there are
$d+1$ points
$\boldsymbol{x}_1,\boldsymbol{x}_2,\ldots,\boldsymbol{x}_{d+1}$ such
that $\boldsymbol{x}_i\in X_i$ for each $i$ and such that
$\mbox{\textup{mpconv}}\{\boldsymbol{x}_1,\boldsymbol{x}_2,\ldots,\boldsymbol{x}_{d+1}\}$
contains the point $\boldsymbol{p}$.
\end{theorem}

\begin{proof}
Write $\boldsymbol{p}=(p_1,\ldots,p_d)$.

There are $\lambda_1,\ldots,\lambda_{d+1}$ and
$\boldsymbol{y}^{(1)},\ldots,\boldsymbol{y}^{(d+1)}\in X_1$ such
that
$$\left\{\begin{array}{l}\max_{j\in[d+1]}\lambda_j=0 \\ \max_{j\in[d+1]}\left(\lambda_j+ \boldsymbol{y}^{(j)}\right)=\boldsymbol{p}.\end{array}\right.$$ Let us read the second equality only for the first component. It means that there is a $j$ such that the first component of $\lambda_j+\boldsymbol{y}^{(j)}$ is equal to $p_1$. Define $\boldsymbol{x}_1$ to be this $\boldsymbol{y}^{(j)}$ and $\mu_1$ to be the corresponding $\lambda_j$. Note that one has $\mu_1+\boldsymbol{x}_1\leq\boldsymbol{p}$ (componentwise), with equality for the first component.

Do the same thing for all the $X_i$ up to $i=d$. Note that one has
then $\mu_i+\boldsymbol{x}_i\leq\boldsymbol{p}$, with equality for
the $i$th component. Now do the same thing for $X_{d+1}$ and define
$\boldsymbol{x}_{d+1}$ to be the $\boldsymbol{y}^{(j)}$ such that
$\lambda_j=0$. Let $\mu_{d+1}:=0$. Note that in this case one has
also $\mu_{d+1}+\boldsymbol{x}_{d+1}\leq\boldsymbol{p}$.

One has then $d+1$ points $\boldsymbol{x}_i\in X_i$ for $i\in[d+1]$
such that $$\left\{\begin{array}{l}\max_{i\in[d+1]}\mu_i=0 \\
\max_{i\in[d+1]}\left(\mu_i+\boldsymbol{x}_i\right)=\boldsymbol{p}.\end{array}\right.$$
\end{proof}


\begin{theorem}\label{thm:trop_gencolcara}
Suppose given $d+1$ finite point sets $X_1, X_2, \ldots, X_{d+1}$
and a max-plus convex set $C$ in $\Rt^d$ such that the max-plus
convex hull of each $X_i$ intersects $C$. Then there are $d+1$
points
$\boldsymbol{x}_1,\boldsymbol{x}_2,\ldots,\boldsymbol{x}_{d+1}$ such
that $\boldsymbol{x}_i\in X_i$ for each $i$ and such that
$\mbox{\textup{mpconv}}\{\boldsymbol{x}_1,\boldsymbol{x}_2,\ldots,\boldsymbol{x}_{d+1}\}$
intersects the max-plus convex $C$.
\end{theorem}

\begin{figure}
\includegraphics[width=10cm]{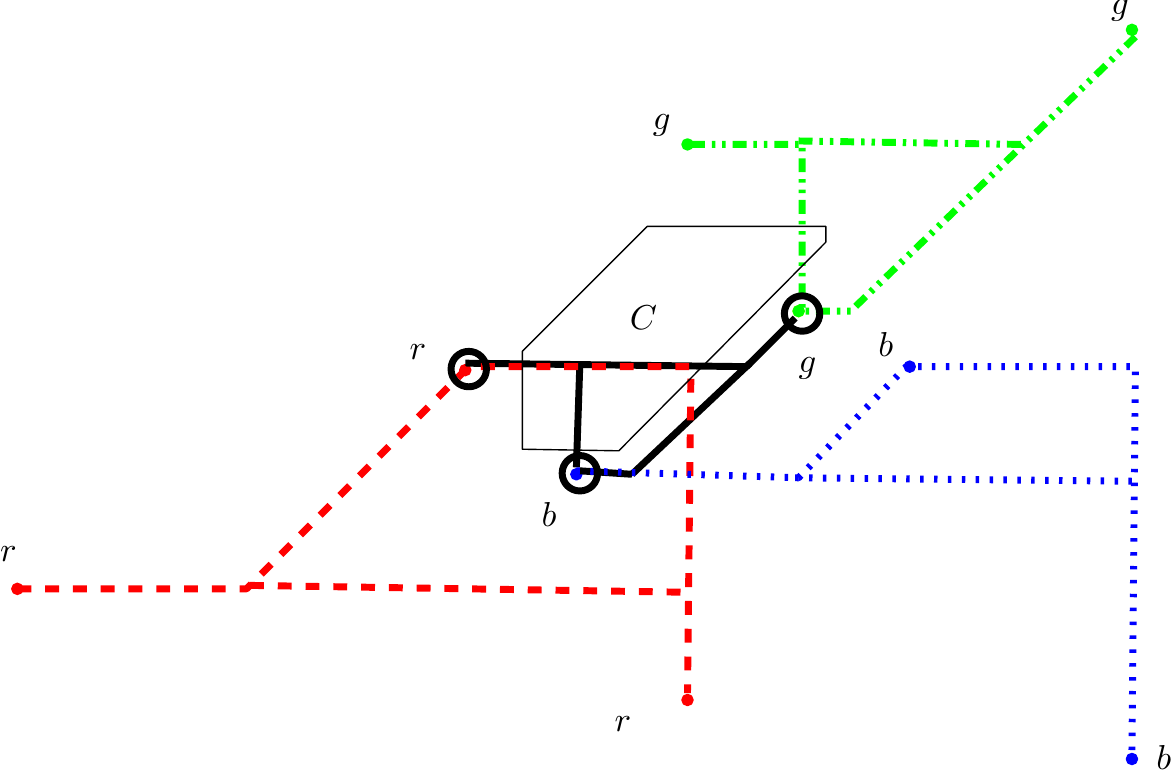}
\label{fig:tropcolorcara} \caption{Illustration of the max-plus
generalized colorful Carathéodory theorem, in dimension 2. There are
three point sets: $X_r$ -- whose points are labelled with $r$ --,
$X_g$ -- whose points are labelled with $g$ -- and $X_b$ -- whose
points are labelled with $b$.}
\end{figure}

Figure \ref{fig:tropcolorcara} is an illustration of this theorem.

To prove this theorem, we will make use of the following lemma:

\begin{lemma}\label{lem:dual}
Consider an $n\times m $ matrix $A=((a_{i,j}))$, with coefficients
in $\Rt$. If $m\geq n$, then, for each column $i$, it is possible to
choose a $\lambda_i\in\Rt$  and add it to each entry of the column
$i$ in such a way that the maxima of each row can be chosen in
different columns. Moreover one can satisfy the additional
requirement $\max_{i=1,\ldots,d+1}\lambda_i=0$.
\end{lemma}

\begin{proof}
The proof works by induction on $n$. If $n=1$, there is nothing to
prove. Hence suppose that $n>1$. Consider the bipartite graph $G$
whose color classes are $W:=[n]$ (the lines of the matrix $A$) and
$U:=[m]$ (the columns of the matrix $A$) and whose edges are those
couples $(i,j)$ such that $a_{i,j}\neq -\infty$. The number
$a_{i,j}$ is then the weight of the corresponding edge. If $G$ has
at least a matching of cardinality $n$, choose such a matching of
maximal weight. The lemma is then a direct consequence of the
duality properties and the complementary slackness for the maximum
assignment problem. Therefore, one can assume that there is no
matching of cardinality $n$. By Hall's marriage theorem, there is a
subset $X$ of $U$ such that $|N(X)|<|X|\leq n$ (where $N(X)$ denotes
the neighborhood of $X$ in $G$). Apply induction on the matrix
having columns $X$ and lines $N(X)$. It provides the values of
$\lambda_i$ on $X$. Define $\lambda_i$ to be $-\infty$ on the other
columns of $A$. One easily checks that one gets the required
property.
\end{proof}

\begin{proof}[Proof of Theorem \ref{thm:trop_gencolcara}] For each $i$, one chooses a point
$\boldsymbol{b}^{(i)}=\left(\begin{array}{c}b_1^{(i)} \\ \vdots \\
b_d^{(i)}\end{array}\right)$ in $C\cap\mbox{mpconv}(X_i)$. Define
$\bar{\boldsymbol{b}}^{(i)}$ to be
$\left(\begin{array}{c}\boldsymbol{b}^{(i)} \\ 0 \end{array}\right)$
and $A$ to be the $(d+1)\times(d+1)$ matrix
$\left(\begin{array}{ccc}\bar{\boldsymbol{b}}^{(1)} & \ldots &
\bar{\boldsymbol{b}}^{(d+1)}\end{array}\right).$

Applying Lemma \ref{lem:dual} on this matrix $A$, one gets that
there is a point $\boldsymbol{p}$ of $C$ such that
$$\left(\begin{array}{c}\boldsymbol{p} \\
0\end{array}\right):=\max_{i=1,\ldots,d+1}\left(\lambda_i+\bar{\boldsymbol{b}}^{(i)}\right),
$$ and such that each component is attained for a different $i$.
Define $\bar{\boldsymbol{p}}:=\left(\begin{array}{c}\boldsymbol{p}
\\ 0\end{array}\right)$.

Now, for each $i$, as $\boldsymbol{b}^{(i)}$ is a max-plus convex
combination of points in $X_i$, one has
$$\bar{\boldsymbol{b}}^{(i)}=\max_{h=1,\ldots,d+1}\left(\mu_h^{(i)}+ \bar{\boldsymbol{a}}_h^{(i)}\right), \qquad\hbox{with}\quad \boldsymbol{a}_h^{(i)}\in X_i\quad\mbox{and}\quad\bar{\boldsymbol{a}}_h^{(i)}:=\left(\begin{array}{c}\boldsymbol{a}_h^{(i)} \\ 0\end{array}\right)\quad\mbox{for all }i,h.$$

There is a $i$ such that the first component of
$\lambda_i+\bar{\boldsymbol{b}}^{(i)}$ is equal to the first one of
$\bar{\boldsymbol{p}}$. One has moreover
$\lambda_i+\bar{\boldsymbol{b}}^{(i)}\leq\bar{\boldsymbol{p}}$
(componentwise). Next, for this $i$, there is a $h(i)\in[d+1]$ such
that
$\mu_{h(i)}^{(i)}+\bar{\boldsymbol{a}}_{h(i)}^{(i)}\leq\bar{\boldsymbol{b}}^{(i)}$
with equality on the first component. Hence, one has
$\lambda_i+\mu_{h(i)}^{(i)}+\bar{\boldsymbol{a}}_{h(i)}^{(i)}\leq\bar{\boldsymbol{p}}$
with equality on the first component.

Do the same thing for all component of $\bar{\boldsymbol{p}}$. Each
component is attained for a different $i$, as noted a few lines
above. Hence one has
$$\max_{i=1,\ldots,d+1}\left(\lambda_i+\mu_{h(i)}^{(i)}+\bar{\boldsymbol{a}}_{h(i)}^{(i)}\right)=\bar{\boldsymbol{p}}.$$
It can be rewritten
$$\boldsymbol{p}=\max_{i=1,\ldots,d+1}\left(\lambda_i+\mu_{h(i)}^{(i)}+\boldsymbol{a}_{h(i)}^{(i)}\right)$$
where $\boldsymbol{a}_{h(i)}^{(i)}$ is a point of $X_i$ for each
$i\in[d+1]$ and
$\max_{i=1,\ldots,d+1}\left(\lambda_i+\mu_{h(i)}^{(i)}\right)=0$.
Define then $\boldsymbol{x}_i:=\boldsymbol{a}_{h(i)}^{(i)}$. A point
of $C$, namely the point $\boldsymbol{p}$, is in their max-plus
convex hull, as required.
\end{proof}

\section{Radon's theorem}

\label{sec:Radon}

The max-plus Radon theorem is :

\begin{theorem}[Max-plus Radon's theorem]\label{thm:mpradon}
Let $X$ be a set of $d+2$ points in $\Rt^d$. Then there are two
pairwise disjoint subsets $X_1$ and $X_2$ of $X$ whose max-plus
convex hulls have a common point.
\end{theorem}

\begin{figure}
\includegraphics[width=10cm]{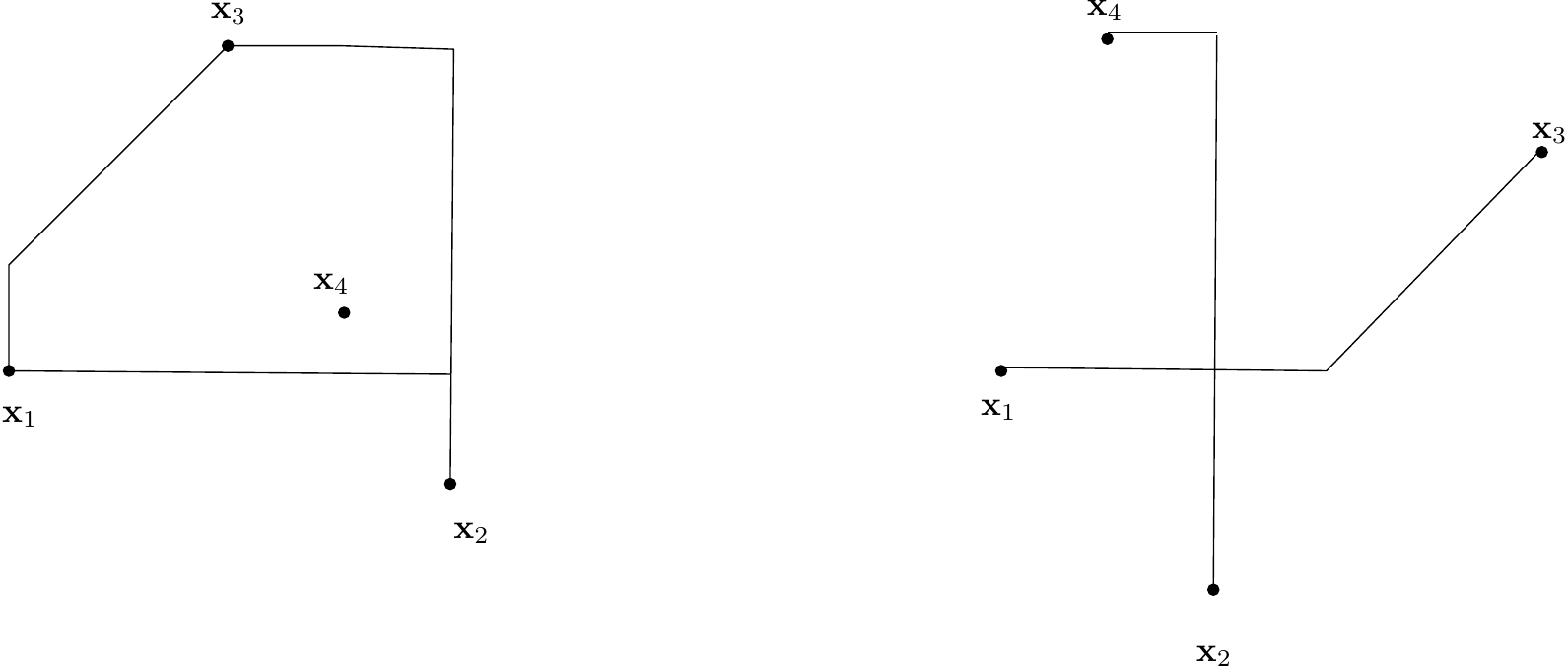}
\label{fig:tropradon}
\caption{Illustration of the max-plus Radon theorem in dimension 2.}
\end{figure}

An illustration is given in Figure \ref{fig:tropradon}.

There are many ways to prove it. When one takes the points in $\R^d$
and not in $\Rt^d$ (no component is equal to $-\infty$), then the
max-plus Radon theorem was proved by Butkovi\v{c} (Theorem 4.7 in
\cite{Butkovic03}) as a consequence of Gondran-Minoux theorem about
max-plus regular matrices \cite{GondranMinoux78} and by Briec and
Horvath in their work about $\mathbb{B}$-convexity
\cite{BriecHorvath04}.

\begin{proof}[Proof of Theorem \ref{thm:mpradon}] We prove the conic version of this theorem:
{\em let $X:=\{\boldsymbol{a}_1,\ldots,\boldsymbol{a}_{d+1}\}$ be a
set of $d+1$ points in $\Rt^d$. Then there are $2$ pairwise disjoint
subsets $S$ and $T$ of $[d+1]$ and $d+1$ reals
$\lambda_1,\ldots,\lambda_{d+1}$ such that $$\max_{i\in
S}\left(\lambda_i+\boldsymbol{a}_i\right)=\max_{i\in
T}\left(\lambda_i+\boldsymbol{a}_i\right).$$ }
The convex version is
derived by adding a 0 component to each point (one gets $d+2$ points
in $\Rt^{d+1}$ and one applies the conic version for the dimension
$d+1$).

\smallskip

Note $\boldsymbol{a}_i=\left(\begin{array}{c} a_{1,1} \\ \vdots \\
a_{d,1} \end{array}\right)$ for $i=1,\ldots,d+1$. Consider the
following matrix $$A(u):=\left(\begin{array}{ccc}u^{a_{1,1}} &
\ldots & u^{a_{1,d+1}} \\ \vdots && \vdots \\ u^{a_{d,1}} & \ldots &
u^{a_{d,d+1}}\end{array}\right)$$ where $u^{-\infty}:=0$. This
matrix has square submatrices, whose determinants are functions of
the undetermined $u$. Let $r$ be the size of the biggest square
submatrix whose determinant is nonzero (as a function of $u$). Take
any $r\times (r+1)$ submatrix of $A(u)$ containing this submatrix
and call it $B(u)$. Without loss of generality, one can assume that
the $r+1$ selected columns are the first $r+1$ of $A(u)$. Denote
$B_i(u)$ the square submatrix of $B(u)$ obtained by deleting the
$i$th column.

\smallskip

We claim that
\begin{equation}\label{eq:radon}\sum_{i=1}^{r+1}(-1)^{i-1}\det\left(B_i(u)\right)\left(\begin{array}{c}u^{a_{1,i}}
\\ \vdots \\
u^{a_{d,i}}\end{array}\right)=\boldsymbol{0}.\end{equation} Indeed,
it is enough to check the validity of this relation for the $j$th
component, for $j=1,\ldots,d$: take the $j$th line of $A(u)$, which
is the line $(u^{a_{j,1}},\ldots,u^{a_{j,d+1}})$, keep the first
$r+1$ components and ``add'' it on the top of $B(u)$ in order to get
a $(r+1)\times(r+1)$ matrix, whose determinant is 0 (still as a
function of $u$):
$$\det\left(\begin{array}{ccc}u^{a_{j,1}} & \ldots & u^{a_{j,r+1}} \\ & B(u) & \end{array}\right)=0.$$
Equation (\ref{eq:radon}) is the development of this equality according to its first line, for $j=1,\ldots,d$.
Remark that the $\det\left(B_i(u)\right)$ are not all equal to 0, by definition of $B(u)$.

\smallskip

Hence, there are $d+1$ functions $p_i(u)=\sum_k \alpha_k
u^{\beta_k}$ (with a finite number of nonzero terms),
$i=1,\ldots,d+1$, not all equal to 0, such that
\begin{equation}\label{eq:radon2}\sum_{i=1}^{d+1}(-1)^{i-1}p_i(u)\left(\begin{array}{c}u^{a_{1,i}}
\\ \vdots \\
u^{a_{d,i}}\end{array}\right)=\boldsymbol{0}.\end{equation} For each
$p_i(u)$, take the largest $\beta_k$ such that $\alpha_k$ is nonzero
and denote it by $\lambda_i$. Let $S$ be the set of $i$ such that
this $\alpha_k$ is strictly positive and $T$ be the set of $i$ such
that this $\alpha_k$ is strictly negative. Reading Equation
(\ref{eq:radon2}) only for the maximum powers leads to the equality
$$\max_{i\in S}\left(\lambda_i+\boldsymbol{a}_i\right)=\max_{i\in T}\left(\lambda_i+\boldsymbol{a}_i\right).$$
\end{proof}

\subsection*{Remark:} In \cite{AkBaGa06}, p. 25-13, Akian, Bapat and Gaubert give various definitions for the rank of a max-plus matrix, and inequalities between them. One gets an one-line proof of the conic max-plus Radon theorem when one combined two of these inequalities:
$$\mbox{rk}_{\mbox{\tiny{GMc}}}(X)\leq\mbox{rk}_{\mbox{\tiny{Schein}}}(X)\leq\mbox{rk}_{\mbox{\tiny{row}}}(X)\leq d.$$

Indeed, $\mbox{rk}_{\mbox{\tiny{GMc}}}(X)$ is the maximum number of
columns that are {\em independent} in the Gondran-Minoux sense
\cite{GondranMinoux02}, that is the maximum number of columns such
that the intersection of the conic hulls of any pairwise disjoint
subsets $S$ and $T$ of these columns is trivial.

\section{Helly's theorem}

\label{sec:Helly}

Radon's theorem can be used to prove a classical theorem about
convex sets, namely Helly's theorem. In the max-plus settings, it is
also possible to derive a max-plus Helly theorem from the max-plus
Radon theorem. The max-plus Helly theorem was first proved by
Gaubert and Sergeev, as a consequence of their work on cyclic
projector in max-plus convexity (\cite{GaubertSergeev07}). Without
$-\infty$-components, a max-plus Helly theorem was proved by Briec
and Horvath in \cite{BriecHorvath04}.

The max-plus version is simply

\begin{theorem}[Max-plus Helly's theorem]
Let $\mathcal{F}$ be a finite collection of max-plus convex sets in
$\Rt^d$. If every $d+1$ members of $\mathcal{F}$ have a nonempty
intersection, then the whole collection have a nonempty
intersection.
\end{theorem}

\begin{proof}
Let $C_1,\ldots,C_n$ be $n$ max-plus convex sets in $\Rt^d$ and
suppose that whenever $d+1$ sets among them are selected, they have
a nonempty intersection. The proof works by induction on $n$ and
first assume that $n=d+2$. Define $\boldsymbol{x}_i$ to be a point
in $\cap_{j=1,\,j\neq i}^{d+2}C_j$. One has then $d+2$ points
$\boldsymbol{x}_1,\ldots,\boldsymbol{x}_{d+2}$. If two of them are
equal, then this point is in the whole intersection. Hence, one can
assume that all the $\boldsymbol{x}_i$ are different. By Radon's
theorem, one has two disjoint subsets $S$ and $T$ partitioning
$[d+2]$ such that there is a point $\boldsymbol{x}$ in
$\mbox{mpconv}\left(\cup_{i\in
S}\boldsymbol{x}_i\right)\cap\mbox{mpconv}\left(\cup_{i\in
T}\boldsymbol{x}_i\right)$. This point $\boldsymbol{x}$ is in every
$C_i$.

Indeed, take $j\in [d+2]$. $j$ is either in $S$ or in $T$. Suppose
w.l.o.g. that $j$ is in $S$. By convexity,
$\mbox{mpconv}\left(\cup_{i\in T}\boldsymbol{x}_i\right)$ is then
included in $C_j$, and hence $\boldsymbol{x}\in C_j$. The case
$n=d+2$ is proved.

Suppose now that $n>d+2$ and that the theorem is proved up to $n-1$.
Define $C_{n-1}':=C_{n-1}\cap C_n$. When $d+2$ max-plus convex sets
$C_i$ are selected, they have a nonempty intersection, according to
what we have just proved. Hence, every $d+1$ members of the
collection $C_1,\ldots,C_{n-2},C_{n-1}'$ have a nonempty
intersection. By induction, the whole collection has a nonempty
intersection.
\end{proof}

\section{Tverberg's theorem}

\label{sec:Tverberg}

One has a Tverberg theorem in the max-plus framework:

\begin{theorem}[Max-plus Tverberg's theorem]\label{thm:trop_tver}
Let $X$ be a set of $(d+1)(q-1)+1$ points in $\Rt^d$. Then there are
$q$ pairwise disjoint subsets $X_1, X_2,\ldots,X_q$ of $X$ whose
max-plus convex hulls have a common point.
\end{theorem}

\begin{figure}
\includegraphics[height=8cm]{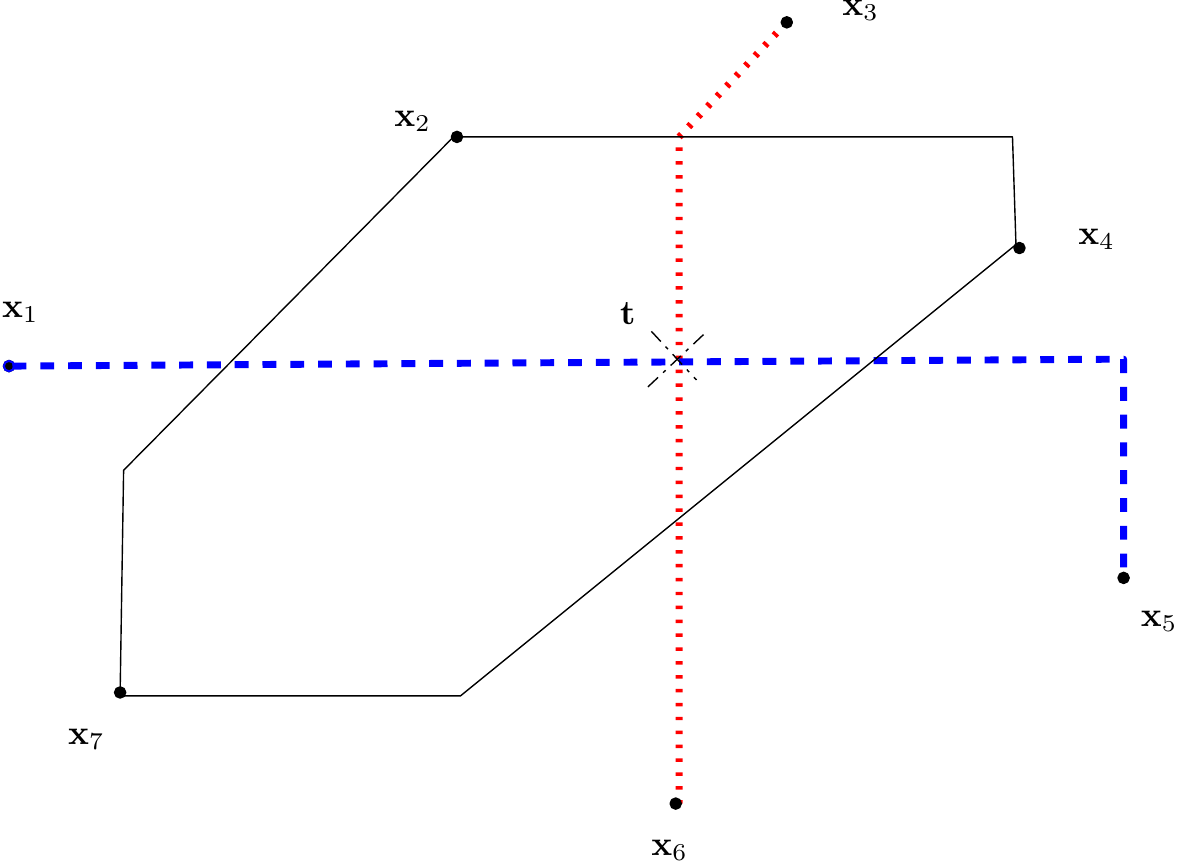}
\label{fig:troptver} \caption{Max-plus Tverberg theorem for $d=2$
and $q=3$.}
\end{figure}

Figure \ref{fig:troptver} illustrates this theorem for $d=2$, $q=3$.
The partition emphasized is
$X_1=\{\boldsymbol{x}_1,\boldsymbol{x}_5\}$,
$X_2=\{\boldsymbol{x}_2,\boldsymbol{x}_4,\boldsymbol{x}_7\}$ and
$X_3=\{\boldsymbol{x}_3,\boldsymbol{x}_6\}$.

To prove this theorem, we will combine the technique used in the
proof of the max-plus Radon theorem above --- identification of
maximum powers in finite series ---, with the beautiful ideas
introduced by Sarkaria \cite{Sarkaria92} and streamlined by
B\'ar\'any and Onn \cite{BaranyOnn97} and Matou\v{s}ek
\cite{Matousek02} to prove the (usual non-max-plus) Tverberg
theorem.

\begin{proof} Put $N:=(d+1)(q-1)$. As for Radon's theorem, we prove the conic version. The convex version is then straightforwardly derived. The conic version is: {\em Let $X=\{\boldsymbol{a}_1,\ldots,\boldsymbol{a}_{N+1}\}$ be a set of $N+1$ points in $\Rt^{d+1}$. Then there are $q$ pairwise disjoint subsets $X_1, X_2,\ldots,X_q$ of $X$ whose max-plus conic hulls have a common
point.}

Write $\boldsymbol{a}_i=\left(\begin{array}{c}a_{1,i} \\ \vdots \\ a_{d+1,i} \end{array}\right)$.

Define linear maps $\phi_j:\R^{d+1}\rightarrow\R^{(d+1)(q-1)}$ for $j\in[q]$ via $$\phi_j(\boldsymbol{y})=(\boldsymbol{0},\ldots,\boldsymbol{0},\boldsymbol{y},\boldsymbol{0},\ldots,\boldsymbol{0})\in\left(\R^{d+1}\right)^{q-1},\quad\hbox{for}\quad j<q,\,\quad \boldsymbol{0}\in\R^{d+1}\quad\hbox{and}\quad \boldsymbol{y}\in\R^{d+1},$$ where $\boldsymbol{y}$ is in $j$th position. Moreover, set $\phi_q(\boldsymbol{y})=(-\boldsymbol{y},-\boldsymbol{y},\ldots,-\boldsymbol{y})$ for $\boldsymbol{y}\in\R^{d+1}$.

For any $\boldsymbol{y}\in\R^{d+1}$, one has $\boldsymbol{0}\in\mbox{conv}\left\{\phi_1(\boldsymbol{y}),\ldots,\phi_q(\boldsymbol{y})\right\}$, in particular for every $i=1,\ldots,N+1$ one has $\boldsymbol{0}\in\mbox{conv}\left\{\phi_1(\boldsymbol{\alpha}_i(u)),\ldots,\phi_q(\boldsymbol{\alpha}_i(u))\right\}$ where $\boldsymbol{\alpha}_i(u)=\left(\begin{array}{c}u^{a_{1,i}} \\ \vdots \\ u^{a_{d+1,i}}\end{array}\right)$, for any real $u$ (one sets $u^{-\infty}:=0$).

Suppose first $u$ fixed. One can apply the colorful Carathéodory theorem to the sets of points $\tilde{X}_1,\ldots,\tilde{X}_{N+1}$, where $\tilde{X}_i:=\{\phi_1(\boldsymbol{\alpha}_i(u)),\ldots,\phi_q(\boldsymbol{\alpha}_i(u))\}$. Indeed, one has $\boldsymbol{0}\in\mbox{conv}\left(\tilde{X}_i\right)$ for each $i$ and we are in $\R^N$. One gets that there exists $j_1,j_2,\ldots,j_{N+1}$ in $[q]$ and  non-negative real numbers $\mu_1,\ldots,\mu_{N+1}$ summing up to 1 such that $$\boldsymbol{0}=\sum_{i=1}^{N+1}\mu_i\phi_{j_i}(\boldsymbol{\alpha}_i(u)).$$

The $j_i$ depend on $u$, of course, but since there are only a
finite number of possible choices, one gets that there exists
$j_1,j_2,\ldots,j_{N+1}$ in $[q]$ and functions
$\mu_1(u),\ldots,\mu_{N+1}(u)$ summing up to 1 such that
\begin{equation}\label{eq:tver}\boldsymbol{0}=\sum_{i=1}^{N+1}\mu_i(u)\phi_{j_i}(\boldsymbol{\alpha}_i(u)).\end{equation}

The $\mu_i(u)$ are solutions of a system of linear equations. Hence
there can be chosen of the form
$$\frac{\sum\gamma_ku^{\beta_k}}{b(u)},$$ with the same denominator
$b(u)$, which is the determinant of the largest invertible
subsystem. Define $S_l:=\{i\in[q]:\,j_i=l\mbox{ and }\mu_i(u)\mbox{
is not constant }=0\}$ (the set of indices $i$ such that $j_i=l$ and
$\mu_i\neq 0$).

Using the definition of the $\phi_j$, one can translate Equation
(\ref{eq:tver}):
$$\sum_{i\in S_1}\mu_i(u)\left(\begin{array}{c}u^{a_{1,i}} \\ \vdots \\ u^{a_{d+1,i}}\end{array}\right)=\ldots=\sum_{i\in S_q}\mu_i(u)\left(\begin{array}{c}u^{a_{1,i}} \\ \vdots \\ u^{a_{d+1,i}}\end{array}\right).$$
Define $\lambda_i$ as the largest $\beta_k$ such that $\gamma_k$ is
non-zero in $\mu_i(u)$ (because they sum up to 1, all the $S_i$ are
non-empty). Reading this equality only for the maximum powers leads
to the equality $$\max_{i\in
S_1}\left(\lambda_i+\boldsymbol{a}_i\right)=\ldots=\max_{i\in
S_q}\left(\lambda_i+\boldsymbol{a}_i\right),\qquad\mbox{where the
}S_l,\,l=1,\ldots,q\mbox{ are disjoint subsets of }[q].$$
\end{proof}

\section{Dutch cheese conjecture}

We finish the article with the max-plus version of Sierksma's
conjecture, which turns out to be a theorem.

\begin{theorem}\label{thm:trop_sierk}
Let $q\geq 2$, $d\geq 1$ and put $N=(d+1)(q-1)$. For every $N+1$
points in $\Rt^d$ the number of unordered max-plus Tverberg
partitions is at least $\left((q-1)!\right)^d$.
\end{theorem}

For instance, if $d=2$ and $q=3$, this theorem says that one has at
least 4 partitions. One can check this assertion in the particular
case given in Figure \ref{fig:troptver}. One partition is
emphasized. There must be three others. Indeed,
$$X_1=\{\boldsymbol{x}_1,\boldsymbol{x}_5\},\, X_2=\{\boldsymbol{x}_2,\boldsymbol{x}_6\},\,X_3=\{\boldsymbol{x}_3,\boldsymbol{x}_4,\boldsymbol{x}_7\},$$
$$X_1=\{\boldsymbol{x}_1,\boldsymbol{x}_4\},\, X_2=\{\boldsymbol{x}_3,\boldsymbol{x}_6\},\,X_3=\{\boldsymbol{x}_2,\boldsymbol{x}_5,\boldsymbol{x}_7\},$$
$$X_1=\{\boldsymbol{x}_1,\boldsymbol{x}_4\},\, X_2=\{\boldsymbol{x}_2,\boldsymbol{x}_6\},\,X_3=\{\boldsymbol{x}_3,\boldsymbol{x}_5,\boldsymbol{x}_7\}$$
are three other Tverberg partitions.

To prove this theorem, we will use a purely combinatorial result --
Corollary \ref{cor:tver_bip} below -- concerning a partition of a
color class of a bipartite graph. It has in a sense a ``Tverberg''
nature. To prove this result, it is useful to prove the following
theorem.

For a graph $G=(V,E)$, let us denote by $N(X)$ the neighborhood of
$X$, that is the set of vertices in $V\setminus X$ having at least
one neighbor in $X$.

\begin{theorem}\label{thm:tverbipq}
Let $G$ be a bipartite graph with color classes $U$ and $W$ and no
isolated vertices, and let $q$ be a positive integer. If $|U|\geq
(q-1)|W|+1$, then there are $q$ disjoint subsets $U_1,\ldots,U_q$ of
$U$ such that $N(U_1)=N(U_2)=\ldots=N(U_q)$. Moreover, there are at
least $\left((q-1)!\right)^{|N(U_1)|-1}$ distinct ways of choosing
these $q$ subsets (one does not take the order into account).
\end{theorem}

\begin{corollary}\label{cor:tver_bip}
Let $G=(V,E)$ be a bipartite graph whose color classes are $U$ and
$W$. Suppose that for all $Y\subseteq W$, $Y\neq\emptyset$, one has
$$|N(Y)|\geq (q-1)|Y|+1.\qquad\qquad\mbox{\textup{(*)}}$$ Then $U$
can be partitioned into $q$ subsets $U_1,\ldots,U_q$ such that for
all $i\in[q]$ one has $N(U_i)=W$. Moreover, there are at least
$\left((q-1)!\right)^{|W|-1}$ distinct partitions satisfying this
property (one does not take the order into account).
\end{corollary}

\begin{proof}[Proof of Theorem \ref{thm:tverbipq}] The proof works by induction on $|U|$. If $|U|=q$, the theorem is clearly true. Hence, let $|U|\geq q+1$.
One can assume that for all $X\subsetneq U$, one has $|X|\leq
(q-1)|N(X)|$ (if not apply induction). One will prove that there
exists a partition $U_1,\ldots, U_q\subseteq U$ such that $N(U_i)=W$
for all $i=1,\ldots,q$, and explains why in this case
$\left((q-1)!\right)^{|W|-1}$ is a lower bound of the distinct ways
of choosing these $q$ subsets, when one does not take the order into
account.

Choose a subset $U'\subseteq U$ of cardinality $(q-1)|W|$. One can
apply Hall's marriage theorem and get a subset of edges $F\subseteq
E$ such that $\deg_F(w)=q-1$ for each $w\in W$ and $\deg_F(u)=1$ for
each $u\in U'$ (make $q-1$ copies of each vertex $w$ of $W$ to see
it). Note that there is no subset $A\neq\emptyset$ of $W$ such that $|N(A)|\leq
(q-1)|A|$, otherwise one would have a subset $X:=U\setminus
N(A)\subseteq U$ such that $(q-1)|N(X)|<|X|$ since $N(X)\subseteq
W\setminus A$. Hence, it is possible to find an order
$w_1,\ldots,w_{|W|}$ of the vertices of $W$ such that $$N(w_i)\cap
\left(Y_1\cup\ldots\cup Y_{i-1}\cup (U\setminus
U')\right)\neq\emptyset\quad\mbox{ for all
$i=1,\ldots,|W|$,}\quad\mbox{(**)}$$ where $Y_j$ denotes the
neighbors of $w_j$ in $F$. In the case when $i=1$, one requires
simply that $N(w_1)\cap (U\setminus U')\neq\emptyset$.

We define now the $U_i$ in order to have $N(U_i)=W$ for
$i=1,\ldots,q$ by adding vertices. Start with
$U_1=\ldots=U_q=\emptyset$.

Add $U\setminus U'$ to $U_1$. Put the $q-1$ vertices of $Y_1$
respectively in $U_2,\ldots, U_q$. The vertex $w_1$ is now in the
neighborhood of $U_1,\ldots, U_q$. Process $w_2,\ldots,w_{|W|}$ in
this order. The processing of $w_i$ consists first in finding the
index $j^*$ such that $w_i$ is already in $N(U_{j^*})$. Such a $j^*$
exists because of property (**). Second, it consists in adding to
each of the $U_j$, except for $j=j^*$, one of the $q-1$ vertices of
$Y_i$. This ensures that when the processing of $w_i$ is finished,
$w_i$ is in the neighborhood of $U_1,\ldots, U_q$. Since all
vertices of $W$ are eventually processed, one gets
$N(U_1)=\ldots=N(U_q)=W$.

It is easy to see why the lower bound on the number of ways of
choosing these $q$ subsets is true. Indeed, there is $q!$ ways of
processing the vertex $w_1$: the subset $U\setminus U'$ can be added
to either subset $U_j$ and the vertices of $Y_1$ to the $q-1$ other
subsets $U_j$ in any order. For each of the vertices $w_i$, there
are $(q-1)!$ ways for adding the vertices of $Y_i$ to each of the
remaining $U_j$. One gets $q!\left((q-1)!\right)^{|W|-1}$ different
ways. Since one does not take the order into account, one has the
required lower bound.
\end{proof}

\begin{proof}[Proof of Corollary \ref{cor:tver_bip}]
Apply Theorem \ref{thm:tverbipq}. One gets $q$ disjoint subsets
$U^{(1)}_1,\ldots,U^{(1)}_q$ having the same neighborhood in $W$.
Define $U':=U\setminus\bigcup_{i=1}^q U^{(1)}_i$ and $W':=W\setminus
N(U^{(1)}_1)$. For all $Y\subseteq W'$, one has $N(Y)\subseteq U'$.
Hence one can apply Theorem \ref{thm:tverbipq} on the subgraph
induced by $U'\cup W'$, and get $U^{(2)}_1,\ldots, U^{(2)}_q$
disjoint subsets of $U'$ having the same neighborhood in $W'$. And
so on. At the end, just define $U_i$ to be the union of all
$U^{(j)}_i$ defined through this process, for each $i=1,\ldots,q$.

The lower bound for the number of distinct partitions is easily
derived.
\end{proof}

We will soon prove Theorem \ref{thm:trop_sierk}. In the proof, we
will need to prove that we have a condition that translates into
condition (*) of Corollary \ref{cor:tver_bip}. This is done by the
following lemma:

\begin{lemma}\label{lem:pert}
For a set
$X=\{\boldsymbol{x}_1,\boldsymbol{x}_2,\ldots,\boldsymbol{x}_n\}$ of
generic points in $\Rt^{d+1}$, if $n\leq (d+1)(q-1)$, then there is
no max-plus conic Tverberg partition into $q$ disjoint subsets.
\end{lemma}

\begin{proof} The proof works by contradiction. Suppose that one has a max-plus conic Tverberg partition $X_1,\ldots,X_q$. Define $\boldsymbol{x}$ to be a Tverberg point, that is a point in the common intersections of the convex hulls of the $X_i$. Define moreover $\lambda_i$ to be the coefficient of $\boldsymbol{x}_i$ in the Tverberg partition.

Consider the graph $H=(V,E)$ where $V:=X$ and there is an edge
between $\boldsymbol{x}_i$ and $\boldsymbol{x}_j$ if the following
two conditions are satisfied: (i) $\boldsymbol{x}_i$ and
$\boldsymbol{x}_j$ are in two consecutive subsets, that is there is
an $l\in\{1,\ldots,q\}$ such that $\boldsymbol{x}_i\in X_l$ and
$\boldsymbol{x}_j\in X_{l+1}$, and (ii) $\lambda_i+\boldsymbol{x}_i$
and $\lambda_j+\boldsymbol{x}_j$ coincide in at least one
coordinate. Parallel edges are allowed.

For each coordinate, $H$ gets at least $q-1$ edges. Hence $H$ has at
least $(d+1)(q-1)\geq n$ edges, and thus has at least one cycle $C$.
Without loss of generality, let
$C:=(\boldsymbol{x}_1,\ldots,\boldsymbol{x}_c)$, in this order,
where $c$ is the size of the cycle. One writes
$\boldsymbol{x}_i:=(x_{1,i},\ldots,x_{d+1,i})$, and define $j(i)$
the coordinate such that
$\lambda_i+x_{i,j(i)}=\lambda_{i+1}+x_{i+1,j(i)}$. Summing the left
and right-hand-side of the equality leads to the following equality:
$$\sum_{i\in C}x_{i,j(i)}=\sum_{i\in C}x_{i+1,j(i)}.$$ Since $j(1),j(2),\ldots,j(c)$ are not all equal, otherwise the edges provided by a coordinate would span a cycle (which would contradict condition (i) defining $H$), one of the term on the left-hand-side of the equality does not appear on the right-hand-side. But then the equality is in contradiction with the genericity assumption.
\end{proof}

\begin{proof}[Proof of Theorem \ref{thm:trop_sierk}.]

We work with the conic version, for there is one-to-one
correspondence between the max-plus Tverberg partitions in the conic
and the convex settings.

Let us start with a particular max-plus Tverberg partition, which
exists because of Theorem \ref{thm:trop_tver}. One has
$X_1,\ldots,X_q$ that provide a partition of
$X=\{\boldsymbol{x}_1,\ldots,\boldsymbol{x}_{N+1}\}$. We are in
$\Rt^{d+1}$. Define $\boldsymbol{x}$ to be a Tverberg point, that is
a point in the common intersections of the convex hulls of the
$X_i$. Define moreover $\lambda_i$ to be the coefficient of
$\boldsymbol{x}_i$ in the Tverberg partition.

Consider the following bipartite graph with color classes $U:=X$ and
$W:=[d+1]$. One puts an edge between $\boldsymbol{x}_i\in X$ and
$j\in[d+1]$ if the $j$th components of $\lambda_i+\boldsymbol{x}_i$
and $\boldsymbol{x}$ coincide.
Our max-plus Tverberg partition $X_1,\ldots, X_q$ provides a
partition of $U$ into subsets $U_1,\ldots,U_q$ such that for all $i$
one has $N(U_i)=W$ by putting $U_i:=X_i$. Moreover, by a slight
perturbation, according to Lemma \ref{lem:pert}, one gets that one
needs at least $(d'+1)(q-1)+1$ points to have a conic Tverberg
partition in dimension $d'+1$, for any $0\leq d'\leq d$. Hence, our
bipartite graph satisfies all conditions of Corollary
\ref{cor:tver_bip}.

Now remark that each partition of $U$ into subsets $U_1,\ldots,U_q$
such that $N(U_i)=W$ for each $i=1,\ldots,q$ provides a max-plus
Tverberg partition by putting $X_i:=U_i$. Corollary
\ref{cor:tver_bip} implies thus the required lower bound for the
number of max-plus Tverberg partitions.
\end{proof}

\begin{small}
\bibliographystyle{amsplain}
\bibliography{tropical}
\end{small}

\end{document}